\documentclass[12pt, a4paper]{article}
\usepackage{amsmath}
\usepackage{amsthm}
\usepackage{amsfonts}
\usepackage{amssymb}
\usepackage{latexsym}

\newcommand{\N}{\ensuremath{\mathbb N}}

\renewcommand{\P}{{\mathbb P}}

\newcommand{\A}{axiom ${\cal A}$ }

\renewcommand{\deg}{\partial}

\theoremstyle{definition}
\newtheorem{df}{Def.}[section]

\theoremstyle{plain}
\newtheorem{sa}[df]{Theorem}
\newtheorem{lem}[df]{Lemma}
\newtheorem{prop}[df]{Proposition}
\newtheorem{cor}[df]{Corollary}

\title{{\bf Meissel's theorem in additive arithmetical semigroups}}
\author{{\bf Stefan Wehmeier}\\ {\small Universit\"at Paderborn, Germany}}
\date{}

\begin{document}
\maketitle

\begin{abstract}
We show how to control the error term in Mertens' formula and related theorems in the context of additive arithmetical semigroups and carry over an old related result of Meissel. 


\end{abstract}

\section{Introduction}

Many algebraic and combinatorial structures may be viewed as multisets of indecomposable components. As a general framework, John Knopfmacher defined the notion of \emph{additive arithmetical semigroup} in his influential books \cite{Kn1} and \cite{Kn2}.

\smallskip

An additive arithmetical semigroup is a free commutative monoid $G$ with non-empty set of generators $P$, equipped with a non-negative integer-valued degree mapping $\deg$ such that
\[ \deg(ab) = \deg(a) + \deg(b) \] 
for all $a, b \in G$, and such that 
\[ G(n) := \# \{ a \in G; \deg(a) = n \} < \infty \] 
for all $n \in \N_0$. Equivalently, we may require that
\[ P(n) := \# \{ p \in P; \deg(p) = n \} \]
be finite for all $n$. In particular, $G$ is countable.

\smallskip

We denote the generating function by
\[ Z(y) := \sum_{n=0}^{\infty} G(n)y^n . \]
Its radius of convergence $\rho$ must satisfy $0 \leq \rho \leq 1$ since $\# G = \sum_n G(n)$ is infinite. Any sensible classification must treat the cases $\rho=0$ and $\rho=1$ separately; in this paper, we restrict our attention to the case $0 < \rho < 1$. Following the traditional notation by Knopfmacher, we set
\[ q:= 1/\rho . \]

\medskip

In the past decades, many authors have investigated the connection between $P(n)$ and $G(n)$. 
A usual hypothesis is that 
\begin{equation}
\label{AxiomA}
\frac{G(n)}{q^n} = A + r(n) 
\end{equation}
 with some $A > 0$ and
\[ r(n) \rightarrow 0 \] sufficiently fast; we term this condition \A (with error term $r(n)$). A desired result would be to show that 
\begin{equation}
\label{andef}
\lambda_n := \frac{n P(n)}{q^n} 
\end{equation}
converges to $1$, in analogy to the prime number theorem; however, this does not necessarily hold even if $r(n) \rightarrow 0$ exponentially fast, see \cite{IndManWar}.    
Wen-Bin Zhang has proved that if $\sum_{n=0}^{\infty} \sup_{k \geq n} |r(k)|$ converges, then $\lambda_n$ is bounded.
On the other hand, in a separate paper (\cite{IndWeh}), Indlekofer and I have shown that examples exist where $r_n = O(\ln(n)^2/n)$ and $\lambda_n$ is unbounded.

However, it is possible to show that the Cesaro limit of $\lambda_n$ exists; more precisely, the book by Zhang and Knopfmacher contains the following theorem (\cite{Kn3}, Section 3.3):
\begin{prop}
\label{TheoremErrorTerm}
Suppose that $\sum_{n=0}^{\infty} \sup_{k \geq n} |r(k)|$ converges.
Then, as $n \rightarrow \infty$,
\begin{itemize}
\item
\begin{equation}
\label{Pkkqkestimate}
 \sum_{k=1}^n \lambda_k = n + O(1). 
\end{equation}
\item 
\begin{equation}
\label{Pkqkestimate2}
 \sum_{k=1}^n \frac{P(k)}{q^k} = \ln(n) + C_1 + O(1/n)
\end{equation}
for some constant $C_1$.
\item
\begin{equation}
\label{ZhangMertens}
\prod_{k=1}^n \left(1 - 1/q^{k} \right)^{P(k)} = \frac{C_2}{n} +
O(1/n^2)
\end{equation}
for some constant $C_2$. 
\end{itemize} 
\end{prop}

In \cite{IndWeh}, Indlekofer and I have shown that Zhang's estimates hold under a much weaker condition, at the price of obtaining no error term:

\begin{prop}
\label{IndWehProp}
Suppose that
\begin{equation}
\label{Hycond}
Z(y) = \frac{H(y)}{1 - qy} 
\end{equation}
and $H(y) \rightarrow A>0$ as $y \rightarrow 1/q$.
Then, as $n \rightarrow \infty$,
\begin{itemize}
\item
\begin{equation}
\label{Lambdakqkestimate}
\sum_{\substack{p \in P, k \in \N \\ \deg(p^k) \leq n}} \frac{1}{k q^{k \deg(p)}} = \ln(n) + \gamma + \ln(A) + o(1)
\end{equation}
where $\gamma$ denotes Euler's constant.
\item
\begin{equation}
\label{Pkqkestimate}
\sum_{k=1}^n \frac{P(k)}{k} = \ln(n) + \gamma + \ln(A) - C_M + o(1) 
\end{equation}
where
\[ C_M = \sum_{k=1}^{\infty} P(k) \left(\frac{1}{q^k} - \ln\left(1 - \frac{1}{q^k}\right)\right) .\]
\item
\[ \sum_{k=1}^n \lambda_k = n + o(n)  . \]
\item
\begin{equation}
\label{Mertenseq}
 \prod_{k=1}^n \left(\frac{1}{1- 1/q^{k}}\right)^{P(k)} = A \exp(\gamma) n + o(n) . 
\end{equation}
\end{itemize}
\end{prop}
In particular, this gives us the values $C_1 = \gamma + \ln(A) - C_M$ and $C_2 = \frac{1}{A \exp(\gamma)}$ for the constants in (\ref{Pkqkestimate2}) and (\ref{ZhangMertens}). Since the proof is based on a Tauberian theorem, it is difficult to turn this into a theorem with explicit error terms.

\section{Mertens-type estimates error terms}

In this section, it will become apparent how the trade-off between hypotheses and error terms in Propositions \ref{TheoremErrorTerm} and \ref{IndWehProp} can be controlled.

\begin{lem}
\label{PkkGquotestimate}
Suppose \A. 
Then
\[ \sum_{k \in \N} k P(k) \frac{G(n-k)}{G(n)} = n - C_3 + o(1), \]
where
\begin{eqnarray*}
 C_3 & := &  \sum_{k=1}^{\infty} k P(k) \sum_{j=2}^{\infty} q^{-jk} \\
   &  = & \sum_{k=1}^{\infty} k P(k) \frac{1}{q^k(q^k-1)}.
\end{eqnarray*}
\end{lem}
\begin{proof}
For all $n \in \N$, we have that
\begin{eqnarray*}
n G(n) & = & \sum_{a \in G_n} \deg(a) \\
& = & \sum_{a \in G_n} \sum_{p^j | a} \deg(p) \\
& = & \sum_p \deg(p) \sum_{j=1}^{\infty} G(n-j\deg(p)) \\
& = & \sum_{k=1}^{\infty} P(k) k \sum_{j=1}^{\infty} G(n-jk).  
\end{eqnarray*}

If $n$ is sufficiently large, we may divide both sides by $G(n)$, 
and \A implies
\[ \sum_{k=1}^{\infty} P(k) k \sum_{j=2}^{\infty} \frac{G(n-jk)}{G(n)} \rightarrow C_3 \]
as $n \rightarrow \infty$ (note that the inner sum is actually finite), such that the lemma follows.
\end{proof}

\begin{sa}
\label{ControlError}
Suppose that $\lambda_n$ is bounded.
Let $r(n) = O(f(n))$ for a continuous nonincreasing function $f$ with $\lim_{x \rightarrow \infty} f(x) = 0$; 
let $F(x):= \int_{1}^{x} f(t) \, dt$. Suppose that $\int_1^{\infty} F(x)/x^2 \, dx$ converges. 
Then 
\begin{itemize}
\item (\ref{Pkkqkestimate}) holds with error term $O(F(n))$.
\item (\ref{Pkqkestimate2}) holds with error term $O(\int_n^{\infty} F(x)/x^2 \, dx)$.
\end{itemize}
\end{sa}
\begin{proof}
First note that $x f(x) = O(F(x))$.
In view of the previous theorem, we have to estimate
\begin{equation}
\label{diff}
\sum_{k=1}^n k P(k) \left(\frac{G(n-k)}{G(n)} - q^{-k}\right) = \sum_{k=1}^n \lambda_k \frac{r(n-k)-r(n)}{A+r(n)}.  
\end{equation}
such that our first assertion follows by reordering the sum $\sum_{k} r(n-k) = \sum_k r(k)$ and comparing it to an integral.

The second assertion follows by partial summation:
Let $S(t):= \sum_{k=1}^{\lfloor t \rfloor} \lambda_k $.
Let $s(t)$ defined by $S(t) = t + s(t) $; then $s(t) = O(F(t))$ by the first part. Let 
\[ I:=  \int_{1}^{\infty} \frac{s(t)}{t^2} \, dt.\]    
Then
\begin{eqnarray*}
\sum_{k=1}^n \lambda_k \cdot \frac{1}{k} & = & S(n) \cdot \frac{1}{n} + \int_1^n S(t) \frac{1}{t^2} \, dt \\   
& = & (n+s(n)) \cdot \frac{1}{n} + \int_1^n \frac{1}{t} \, dt + \int_{1}^n
\frac{s(t)}{t^2} \, dt \\
& = & 1  + O(F(n)/n) +  \ln(n) + I + O\left(\int_{n}^{\infty} \frac{F(t)}{t^2} \, dt\right)  
\end{eqnarray*}
Since $-F(n)/n = \int_n^{\infty} f(x)/x \, dx - \int_{n}^{\infty} F(x)/x^2 \, dx$, it follows that the error term is 
$O(\int_{n}^{\infty} F(x)/x^2 \, dx)$.   
(The main terms are, of course, those of (\ref{Pkqkestimate}), such that
$I + 1 = \gamma + \ln(A) - C_M $.) 
\end{proof}

It follows that Proposition \ref{TheoremErrorTerm} can be improved slightly.
\begin{cor}
The assertions of Proposition \ref{TheoremErrorTerm} hold under the weaker hypotheses
$P(n) = O(q^n/n)$, that $r(n) = O(1/n)$, and that $\sum_{n=1}^{\infty} |r(n)| $ converges.
\end{cor}
\begin{proof}
From our assumptions, it is clear that the right hand side of (\ref{diff}) is bounded in this case. 
The second assertion follows by partial summation as before (this time with $F(x)$ replaced by $1$), and the third one by exponentiation of the second and the definition of $C_M$.
\end{proof}

\section{A sum considered by Meissel}

The following, nearly forgotten result predates Mertens'
  1874 paper \cite{Mertens} containing the theorem named after him. See \cite{LindqvistPeetre2} for a proof and a full account of the historical background.
  
\begin{prop}[Meissel 1866]
For $\alpha \rightarrow 0+ $,
\[ \sum_{p \in \P} \frac{1}{p (\ln(p))^{\alpha}} = \frac{1}{\alpha} + C + O(\alpha), \]
where 
\[ C = \lim_{n \rightarrow \infty} \left(\sum_{p \leq n} \frac{1}{p} -
  \ln(\ln(n))\right). \]
\end{prop}

We show that this result, too, carries over to the case of additive arithmetical semigroups:
\begin{sa}
\label{Meisselsa}
Suppose $H(y) \rightarrow A$ for $ y \rightarrow 1/q$.
Let $C_1 :=  \ln(A) + \gamma - C_M $.
Let $\alpha > 0$. Let 
\[ S(x) := \sum_{k=1}^{\lfloor x \rfloor} \frac{P(k)}{q^k} \] 
and $s(x) := S(x) - \ln(x) - C_1$. Let 
\begin{equation}
\label{Jdef}
J(\alpha) := \int_1^{\infty} \frac{s(x)}{x^{\alpha+1}} \, dx
\end{equation}
Then 
\[ \sum_{k=1}^{\infty} \frac{P(k)}{q^k k^{\alpha}} = 1/\alpha + C_1 + \alpha J(\alpha)  \]
\end{sa}
\begin{proof}
By partial summation. 

\begin{eqnarray*}
\sum_{k=1}^n \frac{P(k)}{q^k k^{\alpha}} & = & S(n) \frac{1}{n^{\alpha}} +
\int_{1}^{n} S(x) \frac{\alpha}{x^{\alpha+1}} \, dx  \\  
& = & \frac{\ln(n) + C_1 + s(n)}{n^{\alpha}} + \alpha \int_{1}^{n}
\frac{\ln(x) + C_1 + s(x)}{x^{\alpha+1}} \, dx.
\end{eqnarray*}
Letting $n \rightarrow \infty$, the first summand on the right vanishes, and the summands of the integral are 
\[
\int_{1}^{\infty} \frac{\ln(x)}{x^{\alpha+1}} \, dx = \frac{1}{\alpha^2},
\]
\[ \int_1^{\infty} \frac{C_1}{x^{\alpha+1}} \, dx =  \frac{C_1}{\alpha}, \]
and (\ref{Jdef}), from which our result follows.
\end{proof}

To arrive at the same result as Meissel, we need $J(0)$ to exist. By Theorem~\ref{ControlError}, it suffices 
if $\int_t^{\infty} F(x)/x^2 dx = O(1/\ln(t)^{1+\epsilon})$ holds.
This finally gives us Meissel's theorem for additive arithmetical semigroups: 
\begin{cor}
If \A{} holds with error term $r(n) = O(1/\ln(n)^{2+\epsilon})$ and $P(n) = O(q^n/n)$, then
\[ \sum_{k=1}^{\infty} \frac{P(k)}{q^k k^{\alpha}} = 1/\alpha + C_1 + O(\alpha)  \]
as $\alpha \rightarrow 0$.
\end{cor}

\bibliography{arisemeng}

\small{
Authors's address: Stefan Wehmeier, Fakult\"at EIM, Universit\"at Paderborn, 33095 Paderborn, Germany. EMail: stefanw@math.upb.de.
}

\end{document}